\documentclass[11pt]{article}
\usepackage{fullpage}
\usepackage{url}
\usepackage{fancybox}
\usepackage{fancyhdr}
\usepackage{tabls}
\usepackage{epsf}
\usepackage{appendix}
\usepackage{ragged2e}
\usepackage[top=1in, bottom=1in, left=1in, right=1in]{geometry}
\usepackage{enumitem}
\usepackage{amsmath, amssymb, amsthm, algorithmic, algorithm}
\usepackage{mathrsfs}
\usepackage{tikz}
\setlist[itemize]{leftmargin=*}
\usepackage{caption}
\usepackage{tabularx}

\newcommand{\apolloR}{APOLLO3\textsuperscript{\textregistered}}
\title{An Application of Reduced Basis Methods to Core Computation in \apolloR}
\author{%
  \textbf{Y.~Conjungo Taumhas$^1$, V.~Ehrlacher$^2$, G.~Dusson$^3$, T.~Lelièvre$^{2}$, and F.~Madiot$^1$}\vspace{3pt} \\
  $^1$Université Paris-Saclay, CEA, \\
  Service d'\'Etudes des Réacteurs et de Mathématiques Appliquées, 91191, Gif-sur-Yvette, France\vspace{6pt} \\ 
  $^2$CERMICS, Ecole des Ponts ParisTech,  \\
  6-8 avenue Blaise Pascal, Cité Descartes, 77455 Marne-la-Vallée, Cedex 2, France\vspace{6pt} \\ 
  $^3$Laboratoire de Mathématiques de
Besancon, UMR CNRS 6623,  \\
Universit\'e de Franche-Comt\'e, 16 route de Gray, 25030 Besan\c con Cedex, France\vspace{6pt}  \\
  \url{yonah.conjungotaumhas@cea.fr}, \url{francois.madiot@cea.fr}, \url{tony.lelievre@enpc.fr},\\
  \url{virginie.ehrlacher@enpc.fr},
  \url{genevieve.dusson@math.cnrs.fr}
}

\newcommand{\keywords}[1]{
  \vspace{-1em}
  \begin{raggedright}
    \begin{list}{}{\setlength{\leftmargin}{0.0in}\setlength{\rightmargin}{0.5in}\small}
    \item[] \ignorespaces KEYWORDS: \upshape #1
    \end{list}
  \end{raggedright}
}

\begin{document}
\maketitle
\justify 
\parskip 6pt plus 1 pt minus 1 pt

\begin{abstract}
In the aim of reducing the computational cost of the resolution of parameter-dependent eigenvalue problems, a model order reduction (MOR) procedure is proposed. We focus on the case of non-self-adjoint generalized eigenvalue problems, such as the stationary multigroup neutron diffusion equations. The method lies in an %linear 
approximation of the manifold of solutions
using a Proper Orthogonal Decomposition approach. 
The numerical method is composed of two stages. In the \textit{offline} stage, we build a reduced space which approximates the manifold. In the \textit{online} stage, for any given new set of parameters, we solve a reduced problem on the reduced space within a much smaller computational time than the required time to solve the high-fidelity problem. This method is applied to core computations in the \apolloR code.
\end{abstract}
\vspace{6pt}
\keywords{Model Order Reduction, Reduced Basis method, Eigenvalue problem, Proper Orthogonal Decomposition}

\section{Introduction} 
We are interested in the parameterized neutron transport equation, when it is solved multiple times for different values of the parameters, e.g. in optimization problems, which is often called a multiquery context. Let us focus on the multigroup approximation over an energy range $[E_{\text{min}},E_{\text{max}}] = [E_G,E_{G-1}] \cup \ldots \cup [E_1,E_0]$, where $G$ stands for the given number of neutron energy groups. Given a parameter $\mu$, the steady-state neutron diffusion equation~\cite[Chapter 7]{DuHa76} seeks the multigroup neutron scalar flux $\phi_{\mu} = \left(\phi_{\mu}^1,\ldots,\phi_{\mu}^G \right)$ associated with the multiplication factor $k_{\text{eff},\mu}$ (the largest eigenvalue in modulus) inside the nuclear reactor core $\mathcal{R}$ such that

\begin{equation}\label{eq:diffusion_eq}
\left(\mathbb{L}^g_{\mu}-\mathbb{H}^g_{\mu}\right)\phi_{\mu}^g = \dfrac{1}{k_{\text{eff},\mu}}\mathbb{F}^g_{\mu}\phi_{\mu}^g,\quad \forall g=\{1,\ldots,G\}, \quad \text{ in }\mathcal{R},%\times \mathbb{S}^2 \times [E_{\text{min}},E_{\text{max}}],
\end{equation}
and vacuum boundary conditions on $\partial\mathcal{R}$ where $\mathcal{R}$ is a bounded and open subset of $\mathbb{R}^3$. The advection operator $\mathbb{L}^g_{\mu}$, the scattering operator $\mathbb{H}^g_{\mu}$ and the fission operator $\mathbb{F}^g_{\mu}$ are defined by
\begin{itemize}
\item $\mathbb{L}^g_{\mu}\phi_{\mu}^g = -\text{div}\left(D^g_{\mu}\nabla\phi_{\mu}^g\right) + \Sigma_{t,\mu}^g\phi_{\mu}^g$;

\item $\mathbb{H}^g_{\mu}\phi_{\mu}^g = \displaystyle\sum_{g'=1}^G\mathbb{H}^{g'\rightarrow g}_{\mu}\phi_{\mu}^{g'}$, where $\mathbb{H}^{g'\rightarrow g}_{\mu}\phi_{\mu}^{g'} = \Sigma_{s,0,\mu}^{g'\rightarrow g}\phi_{\mu}^{g'}$;

\item $\mathbb{F}^g_{\mu}\phi_{\mu}^g = \displaystyle\sum_{g'=1}^G\mathbb{F}^{g',g}_{\mu}\phi_{\mu}^{g'}$, where $\mathbb{F}^{g',g}_{\mu}\phi_{\mu}^{g'} = \chi_{\mu}^g(\nu\Sigma_f)_{\mu}^{g'}\phi_{\mu}^{g'}$;
\end{itemize}
where $D^g_{\mu}$, $\Sigma_{t,\mu}^g$, $\chi_{\mu}^g$, $\nu_{\mu}^{g}$ and $\Sigma_{f_{\mu}}^{g}$ are respectively, for the group $g$, the diffusion coefficient, the total cross-section, the total spectrum, the average number of neutrons emitted per fission, the fission cross-section, and $\Sigma_{s,0,\mu}^{g' \rightarrow g}$ is the Legendre moment of order 0 of the scattering cross-section from group $g'$ to group $g$.
We introduce a partition $(\mathcal{R}_m)_{m=1}^{M}$ of the domain $\mathcal{R}$ with $M\in \mathbb{N}^*$ so that for all $1\leq m \leq M$, $\mathcal{R}_m$ is a domain with Lipschitz, piecewise regular boundaries. For $g,g'=1,2$, the coefficients $D^g_{\mu}$, $\Sigma^g_{t,\mu}$, $\Sigma_{s,0,\mu}^{g' \rightarrow g}$, $\chi^g_{\mu}$, $(\nu\Sigma_f)^g_{\mu}$ are assumed to be piecewise regular on each domain $\mathcal{R}_m$ for $1\leq m \leq M$. 

Several reduced-order models have been proposed in this context~\cite{buchan2013pod,german2019reduced,lorenzi2018adjoint}. 
In this work, we propose a reduced basis (RB) approach, see~\cite{quarteroni2015reduced} for a general introduction and~\cite{CHEREZOV2018195,sartori2015reduced} for applications in neutronics. 
{It is shown numerically in~\cite{yonahPhD} that our approach does not exhibit spurious eigenvalues observed in the monolithic ROM defined in~\cite{german2019reduced} (where the authors propose an alternative groupwise reduced-order model), see~\cite[Remark 2.6]{conjungo2022estimators} for a theoretical argument.}
We focus on the development in the project APOLLO3$^{\text{\textregistered}}$~\cite{mosca2023overview}, a shared platform among CEA, FRAMATOME and EDF, which includes different deterministic solvers for
the neutron transport equation. Particularly, we are interested in the MINARET solver~\cite{lautard_minaret} in the diffusion approximation, 
discretized with discontinuous finite elements.

\section{The reduced basis method} 
\label{sec:reduced_basis}

In this section, we present a methodology for the implementation of a reduced basis solver for a non-self-adjoint eigenvalue problem (see e.g.~\cite{boyaval2010reduced,cances2014greedy}) which writes
\begin{align}
    &   \text{Find $(u_{\mu}, k_{\mu}) \in \mathbb{R}^{\mathcal{N}} \times \mathbb{R}$ such that} \nonumber \\
    &   A_{\mu}u_{\mu} = \frac{1}{k_{\mu}} B_{\mu}u_{\mu},\label{eq:hfp_lambda}
\end{align}
where $k_\mu$ is the largest eigenvalue in modulus, where $\mu\in \mathcal{P}$ stands for the parametric dependence of the problem with $\mathcal{P}$ a compact set of $\mathbb{R}^d, d\ge 1$; $A_{\mu}$ is an invertible non-symmetric matrix, namely the discretized diffusion operator, or disappearance matrix; $B_{\mu}$ is a non-negative non-symmetric matrix, namely the discretized fission operator, or production matrix; and $\mathcal{N}$ is the total number of degrees of freedom of the considered high-fidelity discretization. Typically, for multigroup neutron diffusion calculations, if we denote by $\mathcal{N}_{\mathcal{R}}$ the total number of spatial degrees of freedom, we have $\mathcal{N} = G \times \mathcal{N}_{\mathcal{R}}$.\vspace{6pt}\\

Note that the multiplication factor $k_{\mu}$ is also solution to the following adjoint problem
\begin{align}
    &   \text{Find $(u_{\mu}^*, k_{\mu}) \in \mathbb{R}^{\mathcal{N}} \times \mathbb{R}$ such that} \nonumber \\
    &   A_{\mu}^Tu_{\mu}^* = \frac{1}{k_{\mu}} B_{\mu}^Tu_{\mu}^*,\label{eq:hfp_lambda_adjoint}
\end{align}
where $k_\mu$ is the largest eigenvalue in modulus. The goal is to find a linear space of dimension $n << \mathcal{N}$, denoted by $\mathcal{V}_n$, such that any solution in the manifold 
\[
\displaystyle\mathcal{M}=\{ (u_{\mu},k_{\mu}) ; \mu \in \mathcal{P}\},
\]
can be well-approximated in the space $\mathcal{V}_n$. 

\subsection{Offline Stage} 
To build such a reduced space, we use the information contained in a  training space $\mathcal{P}_{\text{train}}=\left\{\mu_1,\ldots,\mu_{n_s}\right \}$ of $n_s$ parameters. 
The classical \textit{a priori} error analysis exhibits an upper bound on the eigenvalue error which depends on the error on the left and right eigenvectors~\cite{babuvska1991eigenvalue,boffi2010finite}. Following this insight, the reduced space $\mathcal{V}_n$ is built such that
\begin{equation}
    \mathcal{V}_n \subset \text{Span}\left(u_{\mu},u_{\mu}^*;\, \mu \in \mathcal{P}_{\text{train}}\right).
\end{equation}

In order to give the best $n$-rank approximation of the manifold $\mathcal{M}$, we first compute a Singular Value Decomposition (SVD) to the so-called \textit{matrix of snapshots} composed of right eigenvectors
\begin{equation}
    S = \left(u_{\mu_1} | \cdots | u_{\mu_{n_s}} \right) \in \mathbb{R}^{\mathcal{N} \times n_s}, \label{eq:snasphot_matrix}
\end{equation}
which writes
\begin{align}
    &   S = U \Sigma Z^T, \\
    &   U = \left(\xi_1|\ldots|\xi_{\mathcal{N}}\right) \in \mathbb{R}^{\mathcal{N} \times \mathcal{N}}, \nonumber \\
    &   \Sigma = \text{diag}\left(\sigma_1,\ldots,\sigma_{\min (n_s,\mathcal{N})}\right), \nonumber \\
    &   Z = \left(\psi_1|\ldots|\psi_{n_s}\right) \in \mathbb{R}^{n_s \times n_s}, \nonumber
\end{align}
where the $\sigma_i$ are the singular values of $S$, sorted in decreasing order, and $U$ and $Z$ are two orthogonal matrices. Then, the reduced space associated with the approximation of the manifold to Problem~\eqref{eq:hfp_lambda} comes from a Proper Orthogonal Decomposition (POD) and it is defined by
\begin{equation}
    V^{\text{right}} = \left(\xi_1|\ldots|\xi_{n_1}\right),
\end{equation}
where $1\le n_1\le n_s$, which minimizes the 2-norm error between each snapshot and its orthogonal projection onto the subset of dimension $n_1$ spanned by the columns of $V^{\text{right}}$. The integer $n_1$ comes from a truncation of the SVD with respect to a given tolerance criterion related to the singular values. %\cfm{"to a given tolerance criterion"?}
We then proceed similarly to approximate the adjoint manifold $\{ (u^*_{\mu},k_{\mu}) ; \mu \in \mathcal{P}\}$. We perform a SVD to the matrix of snapshots $S^*$  of left eigenvectors, and the POD 
gives the adjoint reduced space
\begin{equation}
   V^{\text{left}} = \left(\xi_1^*|\ldots|\xi_{n_2}^*\right),
\end{equation}
where $1\le n_2\le n_s$.
The resulting reduced space $\mathcal{V}_n$ is defined as the sum of spaces $V^{\text{right}}$ and $V^{\text{left}}$, using an orthonormalization procedure to obtain a basis.

\subsection{Online Stage}

\subsubsection{Assembling the reduced problem}

Let $V_n$ be a matrix containing an orthonormal basis of the reduced space $\mathcal{V}_n$ as columns, a Galerkin projection of Problem~\eqref{eq:hfp_lambda} is done so that, for any $\mu \in \mathcal{P}$, the reduced $n \times n$ matrices are defined as
\begin{align}
 &   A_{\mu,n} = V_n^TA_{\mu}V_n,\\
 &   B_{\mu,n} = V_n^TB_{\mu}V_n.
\end{align}
Assembling such matrices is not trivial, since the high-fidelity sparse matrices $A_{\mu}$ and $B_{\mu}$ are not fully assembled. Indeed,
\begin{align}
    &   A_{\mu} = \begin{pmatrix} A_{\mu}^{1,1} & A_{\mu}^{1,2} & \cdots & A_{\mu}^{1,G} \\ A_{\mu}^{2,1} & A_{\mu}^{2,2} & \cdots & A_{\mu}^{2,G} \\ \vdots & \vdots & \ddots & \vdots \\ A_{\mu}^{G,1} & A_{\mu}^{G,2} & \cdots & A_{\mu}^{G,G}\end{pmatrix} %= \begin{pmatrix} \mathbb{L}^1_{\mu}-\mathbb{H}^{1\rightarrow 1}_{\mu} & -\mathbb{H}^{2\rightarrow 1}_{\mu} & \cdots & -\mathbb{H}^{G\rightarrow 1}_{\mu} \\ -\mathbb{H}^{1\rightarrow 2}_{\mu} & \mathbb{L}^2_{\mu}-\mathbb{H}^{2\rightarrow 2}_{\mu} & \cdots & -\mathbb{H}^{G\rightarrow 2}_{\mu} \\ \vdots & \vdots & \ddots & \vdots \\ -\mathbb{H}^{1\rightarrow G}_{\mu} & -\mathbb{H}^{2\rightarrow G}_{\mu} & \cdots & \mathbb{L}^G_{\mu}-\mathbb{H}^{G\rightarrow G}_{\mu}\end{pmatrix},\\[0.5cm]
    ,
    &   B_{\mu} = \begin{pmatrix} B_{\mu}^{1,1} & B_{\mu}^{1,2} & \cdots & B_{\mu}^{1,G} \\ B_{\mu}^{2,1} & B_{\mu}^{2,2} & \cdots & B_{\mu}^{2,G} \\ \vdots & \vdots & \ddots & \vdots \\ B_{\mu}^{G,1} & B_{\mu}^{G,2} & \cdots & B_{\mu}^{G,G}\end{pmatrix}% = \begin{pmatrix} \mathbb{F}^{1,1}_{\mu} & \mathbb{F}^{2,1}_{\mu} & \cdots & \mathbb{F}^{G,1}_{\mu} \\ \mathbb{F}^{1,2}_{\mu} & \mathbb{F}^{2,2}_{\mu} & \cdots & \mathbb{F}^{G,2}_{\mu} \\ \vdots & \vdots & \ddots & \vdots \\ \mathbb{F}^{1,G}_{\mu} & \mathbb{F}^{2,G}_{\mu} & \cdots & \mathbb{F}_{\mu}^{G,G}\end{pmatrix},
\end{align}
where, for $g'\neq g$, the block matrices $A_{\mu}^{g,g'}$ and $B_{\mu}^{g,g'}$ are sparse $\mathcal{N}_{\mathcal{R}} \times \mathcal{N}_{\mathcal{R}}$ matrices, for $g,g'=\{1,\ldots,G\}$, and the diagonal blocks $A_{\mu}^{g,g}$ are directly accessible in memory.
Therefore, if we decompose the reduced matrix $V_n = \left(\xi_1|\ldots|\xi_n\right)$ along its $G$ group components such that $V_n = \begin{pmatrix}
    \xi_1^1 & \cdots & \xi_n^1 \\ \vdots & & \vdots \\ \xi_1^G & \cdots & \xi_n^G 
\end{pmatrix}$, then we have
\begin{align}
&    (A_{\mu,n})_{i,j} := (V_n^TA_{\mu}V_n)_{i,j}
= \displaystyle\sum_{g=1}^G (\xi_i^g)^TA_{\mu}^{g,g}\xi_j^g + \displaystyle\sum_{\substack{g,g'=1 \\ g'\neq g}}^G (\xi_i^g)^TA_{\mu}^{g,g'}\xi_j^{g'},
\\
&   (B_{\mu,n})_{i,j} := (V_n^TB_{\mu}V_n)_{i,j}
= \displaystyle\sum_{g,g'=1}^G (\xi_i^g)^TB_{\mu}^{g,g'}\xi_j^{g'}
.
\end{align}
For a given parameter $\mu \in \mathcal{P}$, the reduced problem is then the following,% then consists in solving the following problem
\begin{align}
    &   \text{Find $(c_{\mu,n}, k_{\mu,n}) \in \mathbb{R}^{n} \times \mathbb{R}$ such that} \nonumber \\
    &   A_{\mu,n}c_{\mu,n} = \frac{1}{k_{\mu,n}} B_{\mu,n}c_{\mu,n} ,\label{eq:reduced_pb}
\end{align}
where $k_{\mu,n}$ is the largest eigenvalue in modulus. We then obtain 
$
    u_{\mu,n} = V_n c_{\mu,n},
$
as the approximated right eigenvector written in the high-fidelity space $\mathbb{R}^{\mathcal{N}}$.
The associated adjoint problem writes,
\begin{align}
    &   \text{Find $(c_{\mu,n}^*, k_{\mu,n}) \in \mathbb{R}^{n} \times \mathbb{R}$ such that} \nonumber \\
    &   A_{\mu,n}^Tc_{\mu,n}^* = \frac{1}{k_{\mu,n}} B_{\mu,n}^Tc_{\mu,n}^* ,\label{eq:reduced_pb_adjoint}
\end{align}
where $k_{\mu,n}$ is the largest eigenvalue in modulus.
Similarly, we obtain 
$
   u_{\mu,n}^* = V_n c_{\mu,n}^*,
$
as the approximated left eigenvector written in the high-fidelity space $\mathbb{R}^{\mathcal{N}}$.
In order to solve the reduced problem, we use a power iteration method with given relative error tolerances and maximum number of iterations.

\subsubsection{Computing errors and error estimates}

In order to quantify the approximation by the reduced basis method of dimension $n \in \mathbb{N}^*$, we first normalize all high-fidelity and reduced multigroup fluxes such that $\|u_{\mu}\|_2 = \|u_{\mu}^*\|_2 = \|u_{\mu,n}\|_2 = \|u_{\mu,n}^*\|_2 = 1$. Then, we define the respective following $\ell^2$-errors on the eigenvectors and $\ell^2$-error on the eigenvalue
\begin{align}
    &   e_{\mu,n}^u := \|u_{\mu}-u_{\mu,n}\|_2 ,\\
    &   e_{\mu,n}^{u^*} := \|u_{\mu}^*-u_{\mu,n}^*\|_2 ,\\
    &   e_{\mu,n}^k := |k_{\mu}-k_{\mu,n}|.
\end{align}
Let us respectively define the residuals on the direct and adjoint flux by
\begin{align}
    &   R_{\mu,n} := (B_{\mu}-k_{\mu,n}A_{\mu})u_{\mu,n},\\
    &   R_{\mu,n}^* := (B_{\mu}^T-k_{\mu,n}A_{\mu}^T)u_{\mu,n}.
\end{align}
In the following, we will consider the error estimates $\|R_{\mu,n}\|_2$, $\|R_{\mu,n}^*\|_2$ and $\eta^k_{\mu,n} := \dfrac{\|R_{\mu,n}\|\|R_{\mu,n}^*\|}{\langle c_{\mu,n}^*,A_{\mu,n}c_{\mu,n}\rangle}$ respectively on the reduced direct flux $u_{\mu}$, the adjoint reduced flux $u_{\mu}^*$, and the reduced multiplication factor $k_{\mu,n}$~\cite{conjungo2022estimators}. We also introduce the \textit{prefactors} $C_n^u$, $C_n^{u^*}$ and $C_n^k$ by
\begin{align}
    &   C_n^u := \max_{\mu \in \mathcal{P}_{\rm pref}}\; \dfrac{e_{\mu,n}^u}{\|R_{\mu,n}\|_2}, \label{eq:prefactor_u}\\
    &   C_n^{u^*} := \max_{\mu \in \mathcal{P}_{\rm pref}}\; \dfrac{e_{\mu,n}^{u^*}}{\|R_{\mu,n}^*\|_2},\label{eq:prefactor_u_star} \\
    &   C_n^k := \max_{\mu \in \mathcal{P}_{\rm pref}}\; \dfrac{e_{\mu,n}^k}{\eta^k_{\mu,n}},\label{eq:prefactor_k}
\end{align}
where $\mathcal{P}_{\rm pref} \subset \mathcal{P}$ such that $\mathcal{P}_{\rm pref} \cap \mathcal{P}_{\rm train} = \emptyset$.

\section{Numerical tests}

\subsection{Convergence analysis of the RB method on benchmark calculations}

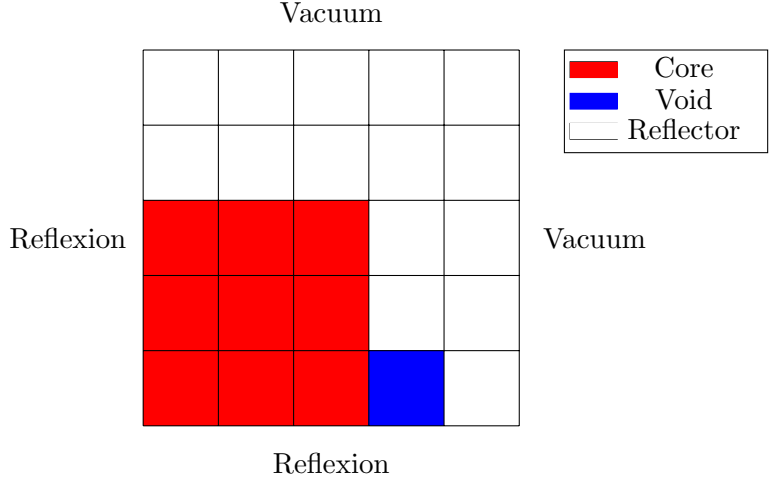
\begin{figure}[htbp]
\begin{center}
\begin{tikzpicture}[scale=5./25.]
%%%%%%%%%%%%%%%%%%%%%%Quarter core%%%%%%%%%%
\fill[white] (0,0) rectangle (25.,25.);
\fill[red] (0,0) rectangle (15.,15.);
\fill[blue] (15.,0.) rectangle (20.,5.);
\draw (0,0) grid[xstep=5.,ystep=5.](25.,25.);
%%%%%%%%%%%%%%%%%%%%%%%%%%%Legend%%%%%%%%%%%%%%%%%%%%%
\draw[black](28.38,23.17) rectangle (31.512,24.2159);
\fill[red] (28.38,23.17) rectangle (31.512,24.2159);
\node at (35.953,23.851925) {Core };
\draw[blue](28.38,21.17-0.1) rectangle (31.512,22.2159-0.1);
\fill[blue] (28.38,21.17-0.1) rectangle (31.512,22.2159-0.1);
\node at (35.953,21.851925-0.1) {Void };
\draw[black](28.38,19.17-0.1) rectangle (31.512,20.2159-0.1);
\fill[white] (28.38,19.17-0.1) rectangle (31.512,20.2159-0.1);
\node at (35.953,19.851925-0.1) {Reflector };
\draw[black](28.,18.17) rectangle (41.512,25.);
%%%%%%%%%%%%%%%%%%%%%%%%%BCS%%%%%%%%%%%%%%%%%%%%%%%%%%%%%%%%%%%
\node at(0.5*25.,-0.1*25.) { Reflexion};
\node at(-0.2*25.,0.5*25.) { Reflexion};
\node at(0.5*25.,+1.1*25.) { Vacuum};
\node at(+1.2*25.,0.5*25.) { Vacuum};
\end{tikzpicture}
\end{center}
\caption{Cross-sectional view of the core ($z = 0$ cm)}
\label{fig:takeda1}
\end{figure}

The reduced basis approach, as implemented in the~\apolloR code, is first tested on Model~1 Case~1 of Takeda neutronics benchmarks~\cite{takeda19913}. A previous work carried out state estimation techniques on this test case using a POD reduced basis from power maps~\cite{conjungocemracs2022}. The considered geometry, as shown in Figure~\ref{fig:takeda1}, is a 3D quarter core in the domain $\{(x,y,z)\in\mathbb{R}^3, 0\leq x\leq 25\,\text{cm}; 0\leq y\leq 25\,\text{cm}; 0\leq z\leq 25\,\text{cm}\} $. The MINARET solver is run with $G=2$ energy groups and $\mathcal{N}_{\mathcal{R}}=3\times 10^5$ spatial degrees of freedom. For this high-fidelity solver, we provide a maximum of 500 outer iterations, with relative $L^2$-error tolerances of $10^{-7}$ and $10^{-8}$ on the two-group flux and on the effective multiplication factor, respectively. The reduced solver runs a power iteration method with respectively relative $\ell^2$-error tolerances of $10^{-8}$ and $10^{-9}$ on the reduced eigenvector and the reduced eigenvalue.

Here, the parameter $\mu$ lies in the 5-dimensional subset $[0.8,1.2]^5$, and then enables small expansions of the equation coefficients such that
\begin{align*}
 &   \left(
D^1_{\mu}, D^2_{\mu},\Sigma_{a,\mu}^1, \Sigma_{a,\mu}^2, \Sigma_{s,0,\mu}^{1\to2}, (\nu\Sigma_f)_{\mu}^1, (\nu\Sigma_{f})_{\mu}^2, \chi^1_{\mu}, \chi^2_{\mu}
\right)\\ = & \left(
\frac{D^1}{\mu_1}, \frac{D^2}{\mu_2},\mu_1\Sigma_{a}^1, \mu_2\Sigma_{a}^2, \mu_3\Sigma_{s,0}^{1\to2}, \mu_4(\nu\Sigma_f)^1, \mu_5(\nu\Sigma_{f})^2, \chi^1, \chi^2
\right),\quad \mu=(\mu_1,\ldots,\mu_5) \in [0.8,1.2]^5,
\end{align*}
where, for $g,g' \in \{1,2\}$, $\Sigma_{a}^g=(\Sigma_{t}^g-\Sigma_{s,0}^{g\to g})$, and the values for the coefficients $D^g$, $\Sigma_{a}^g$, $\Sigma_{s,0}^{g\to g'}$, $(\nu\Sigma_{f})^g$ and $\chi^g$ are given in Appendix 3 of~\cite{takeda19913}.

\begin{figure}[htbp]
\centering
\includegraphics[scale=0.5]{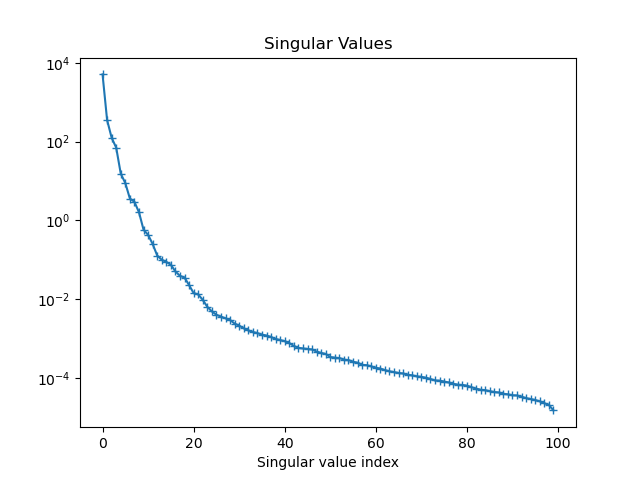}
\includegraphics[scale=0.5]{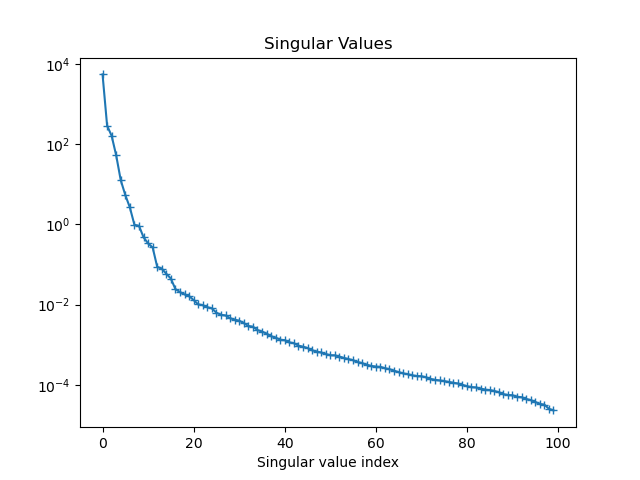}
\caption{Singular values from the SVD of Takeda snapshots, for $n_s=100$. Left: direct eigenvectors; Right: adjoint eigenvectors.} \label{fig:SVD_Takeda}
\end{figure}

We generate a training set $\mathcal{P}_{\text{train}}$ of $n_s=100$ 
parameters with a Latin Hypercube Sampling (LHS) over $[0.8,1.2]^5$. We then compute the SVD of the $2n_s$ \textit{snapshot} matrix as defined in \eqref{eq:snasphot_matrix}. The singular values are shown in Figure~\ref{fig:SVD_Takeda}. The fast decrease of the singular values illustrates the ability of the training set to approximate the manifold of high-fidelity solutions with a reduced basis of small dimension. Here, for example, the 10 first singular values range from $10^4$ to $10^{-1}$. 

\begin{figure}[htbp]
\centering
\includegraphics[scale=0.6]{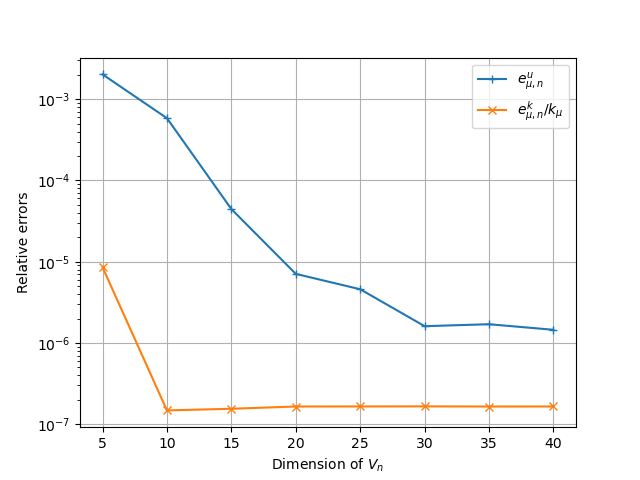}
\caption{Relative errors on the two-group flux and the $k_{\text{eff}}$ with respect to the dimension $n$ of the reduced space $\mathcal{V}_n$, for $\mu=\mu_{\text{test}}$}
\label{fig:Takeda_rel_errors}
\end{figure}

The SVD truncation at the order $n$ then provides a reduced space, and the reduced basis method is tested on the parameter $\mu_{\text{test}}=(1,1,1,1,1)$ which does not belong to the training set %\cgd{correct?} \cyc{yes}
in order to determine to what extend the reduced solver is able to compute a good approximation of the two-group flux and effective multiplicative factor of the Takeda benchmark. The relative errors are depicted in Figure~\ref{fig:Takeda_rel_errors}.
For $n=5$, the reduced solver returns the same $k_{\text{eff}}$ at the order of the pcm ($10^{-5}$), and then the error levels off at the order of magnitude of $10^{-7}$, as the order of convergence is limited by the convergence criterion of the high-fidelity solver.
Note that the test parameter here is particularly well represented by the training space, which explains that the error on the $k_{\text{eff}}$ already reaches the order of the pcm, independently of the order of approximation~$n$. 

\subsection{Computational time reduction on a burnup parametrized nuclear core}

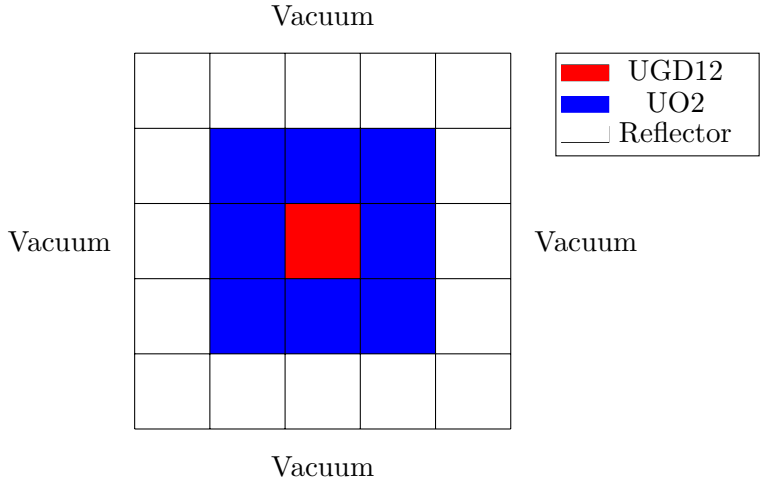
\begin{figure}[htbp]
\begin{center}
%\vspace{-2cm}
\begin{tikzpicture}[scale=5./25.]
%%%%%%%%%%%%%%%%%%%%%%Quarter core%%%%%%%%%%
\fill[white] (0,0) rectangle (25.,25.);
\fill[blue] (5.,5.) rectangle (20.,20.);
\fill[red] (10.,10.) rectangle (15.,15.);
\draw (0,0) grid[xstep=5.,ystep=5.](25.,25.);
%%%%%%%%%%%%%%%%%%%%%%%%%%%Legend%%%%%%%%%%%%%%%%%%%%%
\draw[black](28.38,23.17) rectangle (31.512,24.2159);
\fill[red] (28.38,23.17) rectangle (31.512,24.2159);
\node at (35.953,23.851925) {UGD12 };
\draw[blue](28.38,21.17-0.1) rectangle (31.512,22.2159-0.1);
\fill[blue] (28.38,21.17-0.1) rectangle (31.512,22.2159-0.1);
\node at (35.953,21.851925-0.1) {UO2 };
\draw[black](28.38,19.17-0.1) rectangle (31.512,20.2159-0.1);
\fill[white] (28.38,19.17-0.1) rectangle (31.512,20.2159-0.1);
\node at (35.953,19.851925-0.1) {Reflector };
\draw[black](28.,18.17) rectangle (41.512,25.);
%%%%%%%%%%%%%%%%%%%%%%%%%BCS%%%%%%%%%%%%%%%%%%%%%%%%%%%%%%%%%%%
\node at(0.5*25.,-0.1*25.) { Vacuum};
\node at(-0.2*25.,0.5*25.) { Vacuum};
\node at(0.5*25.,+1.1*25.) { Vacuum};
\node at(+1.2*25.,0.5*25.) { Vacuum};
\end{tikzpicture}
\end{center}
\caption{Median cross-sectional view of the \textit{MiniCore} ($z = 234.36$ cm)}
\label{fig:MiniCore}
\end{figure}

We now test the RB method on a small nuclear core, namely the \textit{MiniCore} problem~\cite{conjungo2022estimators}. The nuclear core geometry is shown in Figure~\ref{fig:MiniCore}. It is a 3D nuclear core and the domain is $\{(x,y,z)\in\mathbb{R}^3, 0\leq x\leq 107.52\,\text{cm}; 0\leq y\leq 107.52\,\text{cm}; 0\leq z\leq 468.72\,\text{cm}\} $. The MINARET solver runs with $G=2$ energy groups and $\mathcal{N}_{\mathcal{R}}=108800$ spatial degrees of freedom. 
For this high-fidelity solver, we provide a maximum of 1000 outer iterations, 
with relative $L^2$-error tolerances of $10^{-7}$ and $10^{-8}$ on the two-group flux and on the effective multiplication factor, respectively. The reduced solver runs a power iteration method with respectively relative $\ell^2$-error tolerances of $10^{-8}$ and $10^{-9}$ on the reduced eigenvector and the reduced eigenvalue.

\begin{figure}[htbp]
\centering
\includegraphics[scale=0.5]{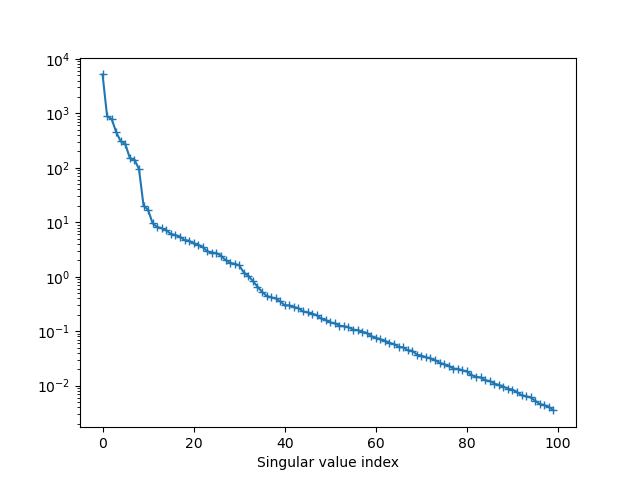}
\includegraphics[scale=0.5]{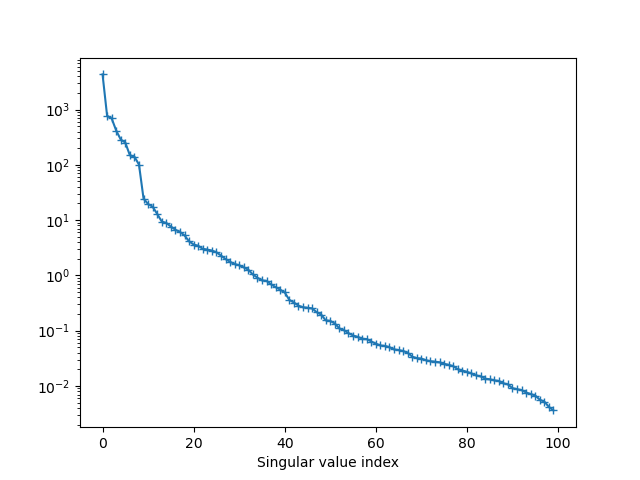}
\caption{Singular values from the SVD of \textit{MiniCore} snapshots, for $n_s=100$. Left: direct eigenvectors; Right: adjoint eigenvectors.} \label{fig:SVD_MiniCore}
\end{figure}

The problem is parametrized by the burnup value for the 9 fuel assemblies (one UGD12 and eight UO2). Here, we generate $n_s=100$ parameters with a Latin Hypercube Sampling (LHS) over the 9-dimensional space
\begin{equation*}
    \mathcal{P}_{\text{train}} \subset \left\{\mu=(\mu_1,\ldots,\mu_9) \in \mathbb{R}^9;\; \mu_1 \in [0,72000];\; \mu_2,\ldots,\mu_9 \in [0,30000]  \right\},
\end{equation*}
where $\mu_1$ is the burnup value of the UGD12 assembly and $\mu_2,\ldots,\mu_9$ are the burnup values of the UO2 assemblies, in MWd/ton. Figure~\ref{fig:SVD_MiniCore} shows an example of SVD with such a training set with $n_s=100$. As in the previous test case, a reduced order model is suitable for the considered snapshot family as the 25 first singular values range from the order of magnitude of $10^4$ to $10^1$.

\begin{figure}[htbp]
\hspace{-2.5cm}
\includegraphics[scale=0.6]{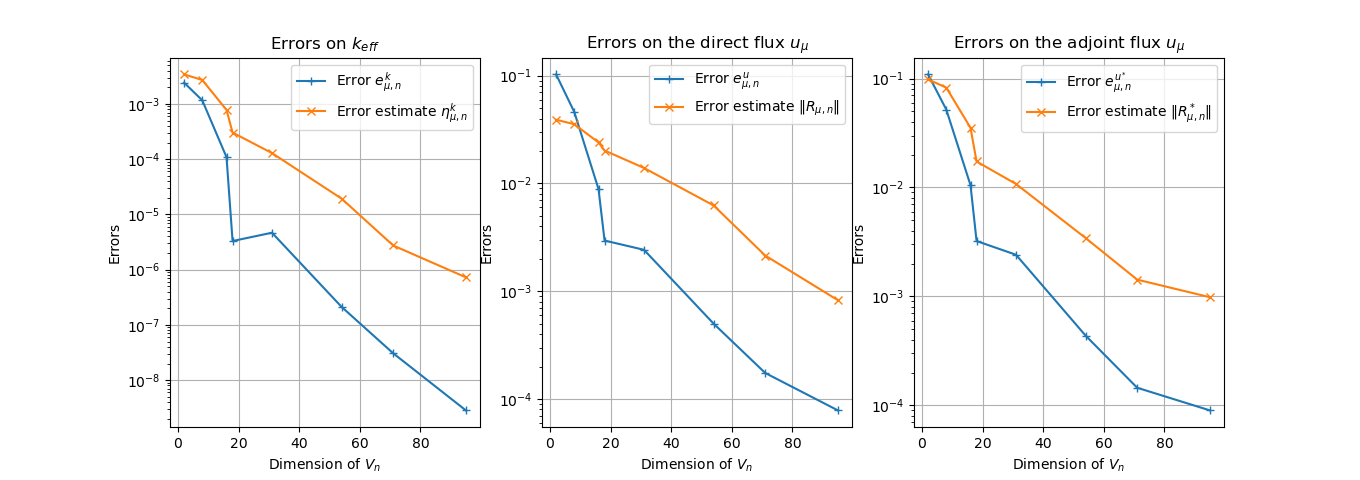}
\caption{Mean errors and their associated error estimates with respect to the dimension $n$ of the reduced space $\mathcal{V}_n$, over $\mathcal{P}_{\text{test}}$. From left to right: error $e_{\mu,n}^k$ and $\eta^k_{\mu,n}$; error $e_{\mu,n}^u$ and residual norm $\|R_{\mu,n}\|_2$; error $e_{\mu,n}^{u^*}$ and residual norm $\|R_{\mu,n}^*\|_2$.}
\label{fig:Errors_MiniCore}
\end{figure}

\begin{figure}[htbp]
\centering
\includegraphics[scale=0.6]{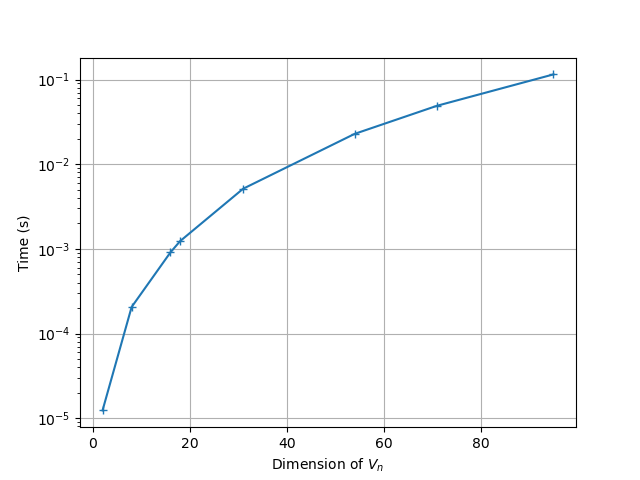}
\caption{Mean computational time for the reduced solver, over $\mathcal{P}_{\text{test}}$}
\label{fig:RB_time}
\end{figure}

\begin{figure}[htbp]
\hspace{-2.5cm}
\includegraphics[scale=0.6]{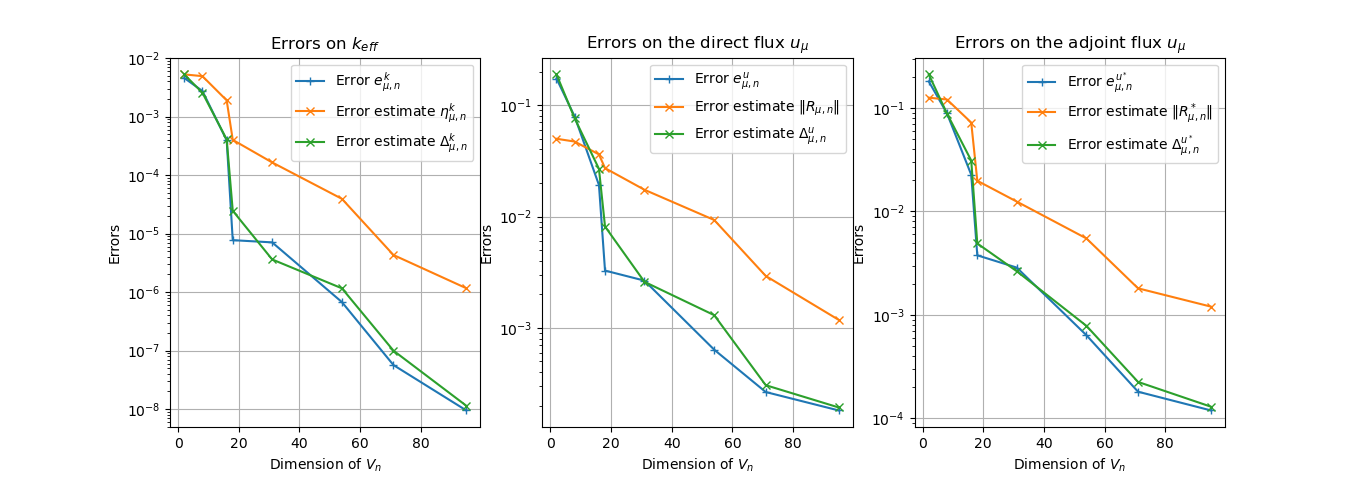}
\caption{Maximum errors and their associated error estimates with respect to the dimension $n$ of the reduced space $\mathcal{V}_n$, over $\mathcal{P}_{\text{test}}$. From left to right: error $e_{\mu,n}^k$, $\eta^k_{\mu,n}$ and $\Delta^k_{\mu,n}$; error $e_{\mu,n}^u$, residual norm $\|R_{\mu,n}\|_2$  and $\Delta^u_{\mu,n}$; error $e_{\mu,n}^{u^*}$, residual norm $\|R_{\mu,n}^*\|_2$ and $\Delta^{u^*}_{\mu,n}$.}
\label{fig:Max_Errors_MiniCore}
\end{figure}

The reduced basis is tested on 10 burnup maps chosen along a LHS over the test space
\begin{equation*}
    \mathcal{P}_{\text{test}} \subset \left\{\mu=(\mu_1,\ldots,\mu_9) \in \mathbb{R}^9;\; \mu_1 = 30000;\; \mu_2,\ldots,\mu_9 \in [0,15000]  \right\}.
\end{equation*}
Figure~\ref{fig:Errors_MiniCore} depicts the convergence of the RB method, as well as the ability for the error estimates to quantify the approximation. We can see that for $n>45$, the errors on the direct and adjoint flux and on the $k_{\text{eff}}$ are respectively below $10^{-3}$ and $10^{-5}$. Regarding the estimates, their convergence, although they are not at the same rate as those of the real errors, are relevant to their potential use in the construction of such an approximation space. Indeed, the cost of the \textit{offline} stage here highly depends on the high-fidelity MINARET solver's cost, as we need to compute high-fidelity solutions for all $\mu \in \mathcal{P}_{\text{train}}$. For the \textit{MiniCore}, one high-fidelity calculation of the $k_{\text{eff}}$ is of the order of the second, whereas computing the reduced eigenvalue $k_{\mu,n}$ requires a computational time of the order of the millisecond, as Figure~\ref{fig:RB_time} shows. In order to get more reliable \textit{a posteriori} error estimates, we define the set $\mathcal{P}_{\rm pref}$ such that
\begin{align*}
    &\mathcal{P}_{\rm pref} \subset \left\{\mu=(\mu_1,\ldots,\mu_9) \in \mathbb{R}^9;\; \mu_1 \in [0,72000];\; \mu_2,\ldots,\mu_9 \in [0,30000]  \right\},\\
    & \#\mathcal{P}_{\rm pref} = 5,\quad \text{and}\quad \mathcal{P}_{\text{pref}} \cap \mathcal{P}_{\text{train}} \cap \mathcal{P}_{\text{test}} = \emptyset .
\end{align*}
We then consider the \textit{a posteriori} error estimates
\begin{align}
    &   \Delta_{\mu,n}^u := C_n^u \|R_{\mu,n}\|_2, \\
    &   \Delta_{\mu,n}^{u^*} := C_n^{u^*} \|R_{\mu,n}^*\|_2, \\
    &   \Delta_{\mu,n}^k := C_n^k \eta^k_{\mu,n},
\end{align}
where the constants $C_n^u$, $C_n^{u^*}$ and $C_n^k$ are defined as in~\eqref{eq:prefactor_u},~\eqref{eq:prefactor_u_star} and~\eqref{eq:prefactor_k} respectively. Figure~\ref{fig:Max_Errors_MiniCore} shows that the estimates defined right below are more reliable as they remain of the same order of magnitude as the real errors, independently of the value of $n$.

\newpage
\section{Conclusion}

The two test cases that were developed highlight the possibility of a reduced basis method implementation in the \apolloR code, in terms of accuracy and computational time reduction. Note that \textit{a posteriori} error estimators in the reduced basis context may be applied in a greedy approach in the \textit{offline} stage~\cite{buffa2012priori,saad2011numerical,giani2016robust}, so that it minimizes calls to the high-fidelity solver, or in an \textit{online} certification of the reduced model~\cite{conjungo2022estimators}. To do so, we should investigate on how to compute the reduced matrices by breaking, as much as possible, their parameter dependency. We could, for example, consider a General Empirical Interpolation Method (GEIM).

\bibliographystyle{siam}
\bibliography{hal_preprint}

\begin{thebibliography}{10}

\bibitem{babuvska1991eigenvalue}
{\sc I.~Babu{\v{s}}ka and J.~Osborn}, {\em Eigenvalue problems}, Handbook of
  Numerical Analysis, 2 (1991), pp.~641--787.

\bibitem{boffi2010finite}
{\sc D.~Boffi}, {\em Finite element approximation of eigenvalue problems}, Acta
  Numerica, 19 (2010), pp.~1--120.

\bibitem{boyaval2010reduced}
{\sc S.~Boyaval, C.~Le~Bris, T.~Leli\`evre, Y.~Maday, N.~C. Nguyen, and A.~T.
  Patera}, {\em Reduced basis techniques for stochastic problems}, Archives of
  Computational methods in Engineering, 17 (2010), pp.~435--454.

\bibitem{buchan2013pod}
{\sc A.~Buchan, C.~Pain, F.~Fang, and I.~Navon}, {\em A pod reduced-order model
  for eigenvalue problems with application to reactor physics}, International
  Journal for Numerical Methods in Engineering, 95 (2013), pp.~1011--1032.

\bibitem{buffa2012priori}
{\sc A.~Buffa, Y.~Maday, A.~T. Patera, C.~Prud’homme, and G.~Turinici}, {\em
  A priori convergence of the greedy algorithm for the parametrized reduced
  basis method}, ESAIM: Mathematical modelling and numerical analysis, 46
  (2012), pp.~595--603.

\bibitem{cances2014greedy}
{\sc E.~Canc{\`e}s, V.~Ehrlacher, and T.~Leli{\`e}vre}, {\em Greedy algorithms
  for high-dimensional eigenvalue problems}, Constructive Approximation, 40
  (2014), pp.~387--423.

\bibitem{CHEREZOV2018195}
{\sc A.~Cherezov, R.~Sanchez, and H.~G. Joo}, {\em A reduced-basis element
  method for pin-by-pin reactor core calculations in diffusion and {SP}3
  approximations}, Annals of Nuclear Energy, 116 (2018), pp.~195--209.

\bibitem{yonahPhD}
{\sc Y.~Conjungo~Taumhas}, {\em Criticality calculations in neutronics: model
  order reduction and a posteriori error estimators}, PhD thesis, Université
  Paris Est, 2023.

\bibitem{conjungo2022estimators}
{\sc Y.~Conjungo~Taumhas, G.~Dusson, V.~Ehrlacher, T.~Leli\`evre, and
  F.~Madiot}, {\em Reduced basis method for non-symmetric eigenvalue problems:
  application to the multigroup neutron diffusion equations}, arXiv e-prints,
  pp.~arXiv--2307.
\newblock https://arxiv.org/abs/2307.05978.

\bibitem{conjungocemracs2022}
{\sc Y.~Conjungo~Taumhas, D.~Labeurthre, F.~Madiot, O.~Mula, and T.~Taddei},
  {\em Impact of physical model error on state estimation for neutronics
  applications}, in ESAIM: Proceedings and Surveys, 2022.

\bibitem{DuHa76}
{\sc J.~J. Duderstadt and L.~J. Hamilton}, {\em Nuclear reactor analysis}, John
  Wiley $\&$ Sons, Inc., 1976.

\bibitem{german2019reduced}
{\sc P.~German and J.~C. Ragusa}, {\em Reduced-order modeling of parameterized
  multi-group diffusion k-eigenvalue problems}, Annals of Nuclear Energy, 134
  (2019), pp.~144--157.

\bibitem{giani2016robust}
{\sc S.~Giani, L.~Grubi{\v{s}}i{\'c}, A.~Miedlar, and J.~S. Ovall}, {\em Robust
  error estimates for approximations of non-self-adjoint eigenvalue problems},
  Numerische Mathematik, 133 (2016), pp.~471--495.

\bibitem{lautard_minaret}
{\sc J.-J. Lautard and J.-Y. Moller}, {\em {MINARET}, a deterministic neutron
  transport solver for nuclear core calculations}, in Proceedings of M\&C 2011,
  Rio de Janeiro (Brazil), May 8-12 2011.

\bibitem{lorenzi2018adjoint}
{\sc S.~Lorenzi}, {\em An adjoint proper orthogonal decomposition method for a
  neutronics reduced order model}, Annals of Nuclear Energy, 114 (2018),
  pp.~245--258.

\bibitem{mosca2023overview}
{\sc P.~Mosca, L.~Bourhrara, A.~Calloo, A.~Gammicchia, F.~Goubioud, L.~Lei-Mao,
  F.~Madiot, F.~Malouch, E.~Masiello, F.~Moreau, S.~Santandrea,
  D.~Sciannandrone, I.~Zmijarevic, E.~Y. Garcias-Cervantes, G.~Valocchi,
  J.~Vidal, F.~Damian, A.~Brighenti, B.~Vezzoni, P.~Laurent, and A.~Willien},
  {\em {APOLLO3}$^{\text{\textregistered}}$: {O}verview of the new code
  capabilities for reactor physics analysis}, in Proceeding of International
  Conference on Mathematics and Computational Methods Applied to Nuclear
  Science and Engineering, M\&C 2023, 2023.

\bibitem{quarteroni2015reduced}
{\sc A.~Quarteroni, A.~Manzoni, and F.~Negri}, {\em Reduced basis methods for
  partial differential equations: an introduction}, vol.~92, Springer, 2015.

\bibitem{saad2011numerical}
{\sc Y.~Saad}, {\em Numerical methods for large eigenvalue problems: revised
  edition}, SIAM: Society for Industrial and Applied Mathematics, 2011.

\bibitem{sartori2015reduced}
{\sc A.~Sartori, A.~Cammi, L.~Luzzi, M.~Ricotti, and G.~Rozza}, {\em Reduced
  order methods: applications to nuclear reactor core spatial dynamics-15566},
  in ICAPP 2015 Proceedings, 2015.

\bibitem{takeda19913}
{\sc T.~Takeda and H.~Ikeda}, {\em 3-{D} neutron transport benchmarks}, Journal
  of Nuclear Science and Technology, 28 (1991), pp.~656--669.

\end{thebibliography}

\end{document}